\documentclass[a4paper, 10pt]{amsart}

\usepackage[a4paper,hmargin=27mm,vmargin=31mm]{geometry}

\usepackage{amsmath,amssymb,hyperref,color}
\hypersetup{hidelinks}

\newtheorem{theorem}{Theorem}[section]
\newtheorem{lemma}[theorem]{Lemma}

\theoremstyle{definition}

\newtheorem{claim}{}[theorem]

\DeclareMathOperator{\Fr}{Frac}
\DeclareMathOperator{\di}{diff}

\renewcommand{\labelenumi}{\arabic{enumi}.}
\numberwithin{equation}{section}

\newcommand{\T}{\mathbb T}
\newcommand{\N}{\mathbb N}

\newcommand{\one}{\mathbf 1}
\DeclareMathOperator{\lcm}{lcm}

\title[Sum-free subsets of sets of integers]{The lower bound for sum-free subsets of sets of integers}
\author{Dongchun Han and Qinghai Zhong}

\address{School of Mathematics, Southwest Jiaotong University, Chengdu, 610031, China}
\email{han-qingfeng@163.com}

\address{Institute of Mathematics and Interdisciplinary Sciences, Xidian University, Xi'an 710126, China}
\email{zhong@xidian.edu.cn} 
\urladdr{https://addcombin.github.io/zhong/}


\begin{document}

\begin{abstract}
We prove that every set of $N\ge2$ positive integers contains a sum-free
subset of cardinality at least $N/3+c\sqrt{\log N}$, where $c>0$ is an
absolute constant. This improves Bedert's lower bound of
$N/3+c_0\log\log N$. 
\end{abstract}
\maketitle

\section{Introduction}\label{s1}

A set $B$ of integers is \emph{sum-free} if $x+y\notin B$ for all
$x,y\in B$. For each positive integer $N$, let $\mathsf S(N)$ be the
largest integer such that every set of $N$ positive integers contains
a sum-free subset of cardinality at least $\mathsf S(N)$.
For background on additive combinatorics and related sum-free set
problems, we refer to the book (\cite{TV}) and survey (\cite{TVS})
of Tao and Vu.
Related extremal problems for finite abelian groups and for
$(k,\ell)$-sum-free subsets were studied by Green and Ruzsa
(\cite{GR}) and Jing and Wu (\cite{JW}), respectively.

In 1965, Erd\H{o}s (\cite[Theorem 2]{Erdos}) proved that
$\mathsf S(N)\ge N/3$ by averaging over dilates of the middle third
of the circle. Alon and Kleitman (\cite[Proposition 1.1]{AK}) subsequently
showed in 1990 that $\mathsf S(N)\ge(N+1)/3$. Bourgain (\cite[Proposition 1.3]{Bourgain}) improved this bound
to $\mathsf S(N)\ge(N+2)/3$ for $N\ge3$. On the other hand, a theorem of
Eberhard, Green and Manners from 2014 implies
$\mathsf S(N)=N/3+o(N)$ (\cite[Theorem 1.1]{EGM}), so the coefficient
$1/3$ is best possible.
The remaining problem is to determine the order of magnitude of the
additional term $\mathsf S(N)-N/3$.

Whether this additional term tends to infinity remained open for many
years. In 2024, Jing and Wu (\cite[Proposition 2.3]{JWJ}) obtained affirmative results
in special cases. The general question
also appears as Problem 1 in the 2025 revision of Green's
\emph{100 open problems} (\cite{Green}). In 2025, Bedert (\cite[Theorem 1.2]{Bedert}) settled the
question by proving that $\mathsf S(N)\ge N/3+c_0\log\log N$ for an
absolute constant $c_0>0$.
In this paper, we improve the additional term from $\log\log N$
to $\sqrt{\log N}$.

\begin{theorem}\label{t1}
There exists an absolute constant $c>0$ such that, for  every set $A$ of $N\ge 2$ positive integers, we have
$\max_{x\in \mathbb R/\mathbb Z}\#\left\{a\in A:ax\in[1/3,2/3)\subset\mathbb R/\mathbb Z\right\}\ge \frac{N}{3}+c\sqrt{\log N}.$ In particular, 
$$\mathsf S(N)\ge N/3+c\sqrt{\log N}.$$
\end{theorem}

\medskip
In Section \ref{s2}, we prove Theorem \ref{t1} assuming the measure
estimate in Lemma \ref{l4}. We partition the set according to the
exponents of $3$ in its elements and combine averaging over three
translates, as in Shakan's argument (\cite[Lemma 2]{Shakan}), with
convexity. To prove Lemma \ref{l4}, we first use estimates of
Koukoulopoulos, Maynard and Yang (\cite{KMY}) in Section \ref{s3} to
count pairs of distinct reduced fractions at small circular distance.
In Section \ref{s4}, we apply this count to the discontinuities of
sums of centered indicator functions and obtain the required measure
estimate.
\section{Proof of the main theorem}\label{s2}

Throughout the manuscript, all logarithms are natural.
For real or complex quantities $X$ and a nonnegative quantity $Y$,
the notation $X\ll Y$, or equivalently $X=O(Y)$, means that
$|X|\le CY$ for some constant $C>0$. Unless otherwise stated,
all implied constants are absolute; in particular, they are
independent of the sets and parameters under consideration.
We write $\#E$ for the cardinality of a finite set $E$.
Let $\N$ be the set of positive integers, let $\N_0=\{0\}\cup \N$, and let $\N_{\ge 2}=\{n\in \N\colon n\ge 2\}$. 

Let $\T=\mathbb R/\mathbb Z$. We shall also use
Markov's inequality (\cite[Section 2.1]{BLM}): for every nonnegative
integrable function $g$ on $\T$ and every real $t>0$, we have
$\mu\{x\in\T:g(x)\ge t\}
\le\frac1t\int_{\T}g\,d\mu$, where $\mu$ is the normalized Lebesgue measure on $\T$.

\newpage

The following lemma, proved in Section \ref{s4}, is a key ingredient
in the proof of Theorem \ref{t1}.

\begin{lemma}\label{l4}
	There exists an absolute constant $\rho>0$ such that every nonempty finite
	set $B\subset\N$ satisfies
	\[
	\mu\left\{x\in\T:\left|\sum_{b\in B}f(bx)\right|\ge\tfrac12\right\}\ge\rho,
	\]
	where  $f=\one_{[1/3, 2/3)}-1/3$ and $\mu$ is the normalized Lebesgue measure on $\T$.
\end{lemma}

We now combine Lemma \ref{l4} with averaging over three translates
and a convexity argument to prove Theorem \ref{t1}. 

\begin{proof}[Proof of Theorem \ref{t1}]
Fix a finite set $A\subset\N$ with $\#A=N$. Let $I=[1/3,2/3)\subset\T$
and $f=\one_I-1/3$.  We define
\[
 T=\max_{x\in\T}\sum_{a\in A}f(ax).
\]
The maximum exists because the sum takes only finitely many values.
For each $x\in\T$,  since
 $I+I=[2/3,1)\cup[0,1/3)\subset \T$ and 
$(I+I)\cap I=\emptyset$,  the set $\{a\in A:ax\in I\}$ is sum-free and
its cardinality is
\[
 \#\{a\in A:ax\in I\}
 =\sum_{a\in A}\one_I(ax)
 =N/3+\sum_{a\in A}f(ax),
\]
whence it suffices to prove 
\begin{center}
	$T\ge c\sqrt{\log N}$ \quad for an absolute
	constant $c>0$.
\end{center}

For each $j\in\N_0$, define
\[
 B_j=\{a\in A:3^j\mid a,\ 3^{j+1}\nmid a\}
 \quad\text{and}\quad
 M_j(x)=T-\sum_{\substack{a\in A\text{ and }3^j\mid a}}f(ax).
\]
For each $u\in\T$, exactly one of the points $u,u+1/3,u+2/3$
lies in $I$, and hence
 $\sum_{r=0}^2f(u+r/3)=0$.
For each $a\in B_j$, write $a=3^jb$ with $3\nmid b$. Since the
residues of $br$ modulo $3$ run through $0,1,2$, we have
\[
 \frac13\sum_{r=0}^2f\bigl(a(x+r/3^{j+1})\bigr)
 =\frac13\sum_{r=0}^2f(ax+br/3)=0.
\]
For each $a\in A$ divisible by $3^{j+1}$, the number $ar/3^{j+1}$
is an integer, so instead
\[
 \frac13\sum_{r=0}^2f\bigl(a(x+r/3^{j+1})\bigr)=f(ax).
\]
Substituting these two identities into the definition of $M_j$, we obtain
\begin{equation}\label{e10}
 \frac13\sum_{r=0}^2M_j\!\left(x+\frac r{3^{j+1}}\right)
 =M_{j+1}(x).
\end{equation}
The definition of $T$ implies $M_0(x)\ge0$. If the function $M_j$
is nonnegative, then its average in \eqref{e10} is also
nonnegative. Induction therefore shows that every function $M_j$
is nonnegative, and hence 
 \begin{equation}\label{e11}
 3M_{j+1}(x)
 =M_j(x)+M_j\!\left(x+\frac1{3^{j+1}}\right)
          +M_j\!\left(x+\frac2{3^{j+1}}\right)\ge M_j(x+\frac{r}{3^{j+1}})\ge 0, \quad  \text{for $r=0,1,2$}.
 \end{equation}

Let $L$ be the number of nonempty sets $B_j$. Choose an integer
$j_0\ge0$ such that $3^{j_0}>\max A$.
Then the function $M_{j_0}$ is identically $T$. For each index $j$
with $B_j=\varnothing$, we have $M_j=M_{j+1}$.
For each of the $L$ indices with $B_j\ne\varnothing$, we have
$M_j(0)\le3M_{j+1}(0)$. Iterating from $j=0$ to $j=j_0-1$, we obtain
\begin{equation}\label{e12}
 T+N/3=M_0(0)\le3^LM_{j_0}(0)=3^LT.
\end{equation}
In particular, $T>0$.
Let $\rho$ be the constant in Lemma \ref{l4}. Then we proceed by the following claim to bound $L$.
\begin{claim}\label{c4}
We have $L\le864\rho^{-3}T^2$.
\end{claim}
\begin{proof}[Proof of \ref{c4}]
Set $K=6T/\rho$ and define the convex function
$\Phi(t)=(\max\{K-t,0\})^2$ for $t\ge0$.
Fix an index $j$ with $B_j\ne\varnothing$.
Let
\[
 E_j=\{x\in\T:M_{j+1}(x)\le K/3\}.
\]
Since each function $x\mapsto f(ax)$ has integral zero, we have
\[
 \int_{\T}M_{j+1}\,d\mu
 =T-\sum_{\substack{a\in A\text{ with }3^{j+1}\mid a}}
       \int_{\T}f(ax)\,d\mu(x)=T.
\]
Applying Markov's inequality to the nonnegative function $M_{j+1}$
with threshold $K/3$, we obtain
\begin{equation}\label{e13}
 \mu(\T\setminus E_j)
 =\mu\{x\in\T:M_{j+1}(x)>K/3\}
 \le\frac3K\int_{\T}M_{j+1}\,d\mu
 =\frac{3T}{K}=\frac\rho2.
\end{equation}

On the other hand, in view of the identity $M_{j+1}-M_j=\sum_{a\in B_j}f(ax)$, it follows from Lemma \ref{l4} (applied to $B_j$)
that
 $\mu\{x\in\T:|M_j(x)-M_{j+1}(x)|\ge1/2\}\ge\rho$,
whence \eqref{e13} implies 
\begin{equation}\label{e14}
 \mu\{x\in E_j:|M_j(x)-M_{j+1}(x)|\ge1/2\}\ge\rho/2.
\end{equation}
By \eqref{e10}, the three nonnegative values $M_j(x+r/3^{j+1})$
have mean $M_{j+1}(x)$. Convexity of the function $\Phi$ therefore implies
\[
 \frac13\sum_{r=0}^2\Phi(M_j(x+r/3^{j+1}))-\Phi(M_{j+1}(x))\ge0.
\]
For each $x\in E_j$, by \eqref{e11}, we have
\[
 0\le M_j(x+r/3^{j+1})\le3M_{j+1}(x)\le K,
 \qquad \text{ for  }r=0,1,2.
\]
All these values, including $M_{j+1}(x)$, lie in $[0,K]$,
where the function $\Phi$ equals $(K-t)^2$. Expanding
the squares and using \eqref{e10}, we obtain
\begin{equation}\label{e15}
 \frac13\sum_{r=0}^2\Phi(M_j(x+r/3^{j+1}))-\Phi(M_{j+1}(x))
 =\frac13\sum_{r=0}^2|M_j(x+r/3^{j+1})-M_{j+1}(x)|^2 \text{ for every }x\in E_j.
\end{equation}
Using translation invariance of $\mu$, discarding the nonnegative
integrand outside $E_j$, and retaining only the term with $r=0$ in
\eqref{e15}, we obtain
\begin{equation}\label{e16}
 \int_{\T}\bigl(\Phi(M_j)-\Phi(M_{j+1})\bigr)\,d\mu
 \ge\frac13\int_{E_j}|M_j-M_{j+1}|^2\,d\mu
 \ge\frac\rho{24}.
\end{equation}
The last inequality follows from \eqref{e14}, since the squared
difference is at least $1/4$ on a set of measure at least $\rho/2$.

For the integer $j_0$ chosen above, the function $M_{j_0}$ is
identically $T$. Among the indices
$0\le j<j_0$, exactly $L$ satisfy $B_j\ne\varnothing$; each contributes
an integral difference of at least $\rho/24$ by \eqref{e16}.
For every other index,  we have $M_j=M_{j+1}$,
so the integral difference is zero. Summing implies that 
\[
 \frac{L\rho}{24}
 \le\sum_{j=0}^{j_0-1}
 \int_{\T}\bigl(\Phi(M_j)-\Phi(M_{j+1})\bigr)\,d\mu=\int_{\T}\Phi(M_0)\,d\mu-\Phi(T)
 \le K^2={36T^2}/{\rho^2}.
\]
The final inequality uses $0\le\Phi(t)\le K^2$ for $t\ge0$ and
$\mu(\T)=1$. Consequently, $L\le864\rho^{-3}T^2$.

\qedhere[End of \ref{c4}]
\end{proof}

Combining \eqref{e12} with \ref{c4}, we obtain
 $\log N\le L\log3+\log(3T)
 \le(864\rho^{-3}\log3+9)T^2$,
where the last inequality uses $\log u\le u^2$ for $u>0$.
This proves the theorem.
\end{proof}

\section{Estimates for pairs of reduced fractions}\label{s3}

We recall the elementary arithmetic functions used below. For  $n\in \N$, \emph{Euler's totient function} and the \emph{divisor function} are
defined by
 $\varphi(n)=\#\{1\le r\le n\colon \gcd(r,n)=1\}$,
 and 
 $\tau(n)=\#\{d\in\N\colon d\mid n\}$.
In particular, $\varphi(1)=\tau(1)=1$. An arithmetic function $v$ is
\emph{multiplicative} if $v(1)=1$ and $v(rs)=v(r)v(s)$ whenever
$\gcd(r,s)=1$. Both $\varphi$ and $\tau$ are multiplicative, with
$\varphi(p^k)=p^{k-1}(p-1)$ and $\tau(p^k)=k+1$ for every prime $p$
and integer $k\ge1$. We write $\mathsf v_p(n)$ for the exponent of $p$ in $n$. We  shall also use the following elementary estimates:
\begin{equation}\label{e1}
 \log\frac n{\varphi(n)}=\sum_{p\mid n}\frac1p+O(1), \quad 
 \sum_{i=1}^n\frac1i\le1+\log n,\quad \text{ and }\ \sum_p\frac1{p(p-1)}
 \le 1.
\end{equation}

For $x,y\in\T$, their circular distance is
 $|x-y|:=\min_{k\in\mathbb Z}|\widetilde x-\widetilde y-k|$,
where $\widetilde x,\widetilde y\in\mathbb R$ are any representatives
of $x,y$. This definition is independent of the representatives and
equals the length of the shorter arc joining $x$ and $y$ on the circle
of circumference one. In particular, $0\le |x-y|\le1/2$.

For a finite subset $\mathcal D\subset \N_{\ge 2}$, let
$\Fr(\mathcal D)\subset\T$ consist of all reduced fractions
with denominator in $\mathcal D$.
For $a,b\in\N$ and every real $t\ge1$, we define
\[
 L_t(a,b):=\sum_{\substack{p\text{ divides } ab/\gcd(a,b)^2  \text{ with } p>t}}\frac1p.
\]

We first state the
two estimates of Koukoulopoulos, Maynard and Yang.

\begin{lemma}\label{l1}
There exists a constant $C_1>0$ such that the following hold.
\begin{enumerate}
\renewcommand{\labelenumi}{\arabic{enumi}.}
\item Let $v$ be a nonnegative multiplicative function satisfying
$v(n)\le\tau(n)$ for every  $n\in \N$. For every real $x\ge1$, we have
\[
 \sum_{h\le x}v(h)
 \le C_1x\exp\!\left(\sum_{p\le x}\frac{v(p)-1}{p}\right),
\]
where the sum in the exponent is over primes.
\item Let $\mathcal D\subset \N_{\ge 2}$ be  nonempty and finite, and let $J=\sum_{a\in\mathcal D}\varphi(a)$.
For all real $y_0,t\ge1$, we have
\[
 \sum_{\substack{a,b\in\mathcal D \text{ with}\\
                 \lcm(a,b)\le y_0\text{ and } L_t(a,b)\ge1}}
 \varphi(a)\varphi(b)
 \le C_1e^{-t}\sqrt {y_0}\,J^{3/2}.
\]
\end{enumerate}
\end{lemma}

\begin{proof}
1. Apply \cite[Lemma 6.1]{KMY} with the ordinary divisor function.

\smallskip
2. The assertion is a direct special case of \cite[Proposition 7.2]{KMY}, obtained by taking
 $\varepsilon=\frac12$, $C=s=1$,
 $Q=\max\mathcal D$, and $\psi(a)=a\,\one_{\mathcal D}(a)$.
For $a\in\mathcal D$, the weight $\psi(a)\varphi(a)/a$ is $\varphi(a)$;
outside $\mathcal D$, it is zero. The sum
$\Psi(Q)$ in that proposition is $2J$, and the quantity $D(a,b)$
defined in \cite[Equation (5.1)]{KMY} equals $\lcm(a,b)$.
The prime sum $L_t(a,b)$ in that proposition agrees with our definition
(see \cite[Equation (5.2)]{KMY}).
This implies the stated estimate, with the factor $2^{3/2}$ absorbed
into the absolute constant.
\end{proof}

We then count pairs of reduced fractions with fixed denominators
whose circular distance is at most $R$.

\begin{lemma}\label{l2}
	There exists a constant $C_0>0$ such that, for any two  $a,b\in \N_{\ge 2}$ and every real $R$ with
	$1/\lcm(a,b)\le R<1/2$, we have
	\[
	\#\left\{\left(\frac{x}{a},\frac{y}{b}\right)\in \T^2\colon
	0<\left|\frac{x}{a}-\frac{y}{b}\right|\le R,\  \gcd(a,x)=\gcd(b,y)=1
	\right\}
	\le C_0 R\varphi(a)\varphi(b)e^{L_{R\lcm(a,b)}(a,b)},
	\]
	where $x,y$ range over integers and
	$\left|\frac{x}{a}-\frac{y}{b}\right|$ denotes circular distance on $\T$.
\end{lemma}

\begin{proof}
	Let $V$ denote the number of ordered pairs in the statement.
	Set
	\[
	a'=a/\gcd(a,b),\quad b'=b/\gcd(a,b),\quad n=a'b',\quad \text{ and }\quad
	m=\lcm(a,b).
	\]
	The difference of two fractions with denominators $a,b$ is a
	multiple of $1/m$. Since $R<1/2$, every pair counted by
	$V$ determines a unique nonzero integer $h$ with $|h|\le mR$ such that
	$\frac xa-\frac yb\equiv\frac hm\pmod1$.
	
	For every nonzero integer $h$, we define
	\[
	c(h):=\#\left\{(x,y)\in
	(\mathbb Z/a\mathbb Z)\times(\mathbb Z/b\mathbb Z)\colon
	b'x-a'y\equiv h\pmod m,\ 
	\gcd(x,a)=\gcd(y,b)=1
	\right\}.
	\]
	Thus, by summing over the possible
	values of $h$, we obtain
	\begin{equation}\label{e2}
	V=\sum_{\substack{h\in\mathbb Z\text{ with }1\le |h|\le mR}}c(h).
	\end{equation}
	
	\medskip
	We proceed by the following claim to bound $c(h)$.
	\begin{claim}\label{c1}
		For every nonzero integer $h$, we have $c(h)=0$ if $\gcd(|h|,n)\neq 1$, and 
		\[
		c(h)\le\frac{\varphi(a)\varphi(b)}{m}\,
		\frac{n}{\varphi(n)}\frac{|h|}{\varphi(|h|)} \quad \text{ if }\gcd(|h|,n)= 1.
		\]
	\end{claim}
	\begin{proof}[Proof of \ref{c1}]
		Let $p$ be a prime divisor of $ab$, let
		$\alpha=\mathsf v_p(a)$, and let $\beta=\mathsf v_p(b)$. Thus $\alpha,\beta\ge0$
		with $\max\{\alpha,\beta\}\ge1$. Define $c_p(h)$ to be the cardinality of the following set: 
		\[
	\left\{(x,y)\in
		(\mathbb Z/p^\alpha\mathbb Z)\times(\mathbb Z/p^\beta\mathbb Z):
		 b'x-a'y\equiv h\pmod{p^{\max\{\alpha,\beta\}}},\ 
		 \gcd(x,p^\alpha)=\gcd(y,p^\beta)=1 \right\}.
		\]
		Since $p^\alpha b'$ and $p^\beta a'$ are both divisible by
		$p^{\max\{\alpha,\beta\}}$,
		the congruence is well-defined on these residue classes.
		It follows from the Chinese Remainder Theorem that
		\begin{equation}\label{e3}
		c(h)=\prod_{p\mid ab}c_p(h).
		\end{equation}
		
		Suppose first that $p\mid n$. Thus $\alpha\neq \beta$, say $\alpha>\beta$ by symmetry, whence $p\mid a'$ and
		$p\nmid b'$. For each $y\pmod {p^\beta}$ with
		$\gcd(y,p^\beta)=1$, the congruence determines
		exactly one $x\pmod {p^\alpha}$, namely
		$x\equiv(b')^{-1}(h+a'y)\pmod{p^\alpha}$, so 
		\begin{center}
			$\gcd(x, p^{\alpha})=1$\quad  if and only if\quad  $p\nmid h$.
		\end{center}
		There are
		$\varphi(p^\beta)$ choices for $y$, including one choice when
		$\beta=0$, since $\varphi(1)=1$. Therefore
		\begin{equation}\label{e4}
		c_p(h)=
		\begin{cases}
			\varphi(p^\beta),&\text{if }p\nmid h,\\
			0,&\text{if }p\mid h.
		\end{cases}
		\end{equation}

		Suppose next that $p\nmid n$. Then $\alpha=\beta=k\ge1$. Both $a'$ and $b'$
		are invertible modulo $p^k$. Each $y\pmod{ p^k}$ with $p\nmid y$ determines
		a unique $x \pmod {p^k}$, namely
		$x\equiv(b')^{-1}(h+a'y)\pmod{p^k}$, so 
		\[
		p\nmid x\quad\text{ if and only if }\quad h+a'y\not\equiv0\pmod p.
		\]
        If $p\mid h$, every such $y$ is admissible. Otherwise, exactly
        one of the $p-1$ nonzero classes of $y\pmod p$ is excluded,
        removing $p^{k-1}$ choices. Therefore
 		\begin{equation}\label{e5}
		c_p(h)=
		\begin{cases}
			p^{k-1}(p-1),&\text{if }p\mid h,\\
			p^{k-1}(p-2),&\text{if }p\nmid h.
		\end{cases}
		\end{equation}
		
		\medskip
		
		We now combine the local counts.
		 If $\gcd(|h|,n)>1$, choose
		a prime $p\mid\gcd(|h|,n)$. It follows from \eqref{e4} that $c_p(h)=0$, whence  $c(h)=0$ by \eqref{e3}.
		
		Assume that $\gcd(|h|,n)=1$. If $p\mid n$, then $p\nmid h$; if $p\mid h$, then $p\nmid n$.  It follows from \eqref{e4} and \eqref{e5} that 
		\[
		\frac{c_p(h)}{\varphi(p^{\mathsf v_p(a)})\varphi(p^{\mathsf v_p(b)})/p^{\max\{\mathsf v_p(a),\mathsf v_p(b)\}}}
		=\begin{cases}
			(1-1/p)^{-1},&\text{if }p\mid n\text{ or }p\mid h,\\
			\dfrac{1-2/p}{(1-1/p)^2}\le1,
			&\text{if }p\nmid n\text{ and }p\nmid h.
		\end{cases}
		\]
		Multiplying over $p\mid ab$ and enlarging the product over primes
		dividing both $ab$ and $h$ to all primes dividing $h$, whose added
		factors $(1-1/p)^{-1}$ are greater than one, we obtain
		\[
		c(h)=\prod_{p\mid ab}c_p(h)
		\le\frac{\varphi(a)\varphi(b)}m
		\prod_{p\mid n}(1-1/p)^{-1}
		\prod_{p\mid h}(1-1/p)^{-1}=\frac{\varphi(a)\varphi(b)}m
		\frac n{\varphi(n)}\frac{|h|}{\varphi(|h|)},
		\]
		as required.
		\qedhere[End of \ref{c1}]
	\end{proof}

\medskip	
Since $mR\ge1$, we may apply Lemma \ref{l1}.1 (with $x=mR$) to a suitable
multiplicative function. Define
\[
 v(h)=\frac h{\varphi(h)}\one_{\gcd(h,n)=1},
 \qquad \text{ where }h\in\N,
\]
which is nonnegative and multiplicative. For every prime $p$
and $k\in \N$, we have $0\le v(p^k)\le p/(p-1)\le2\le k+1$, whence $v(h)\le\tau(h)$ for every $h\in\N$.
For primes $p\nmid n$, we have $(v(p)-1)/p=1/(p(p-1))$,
whose sum is at most one by \eqref{e1}; for primes $p\mid n$, we have
$v(p)=0$. Thus Lemma \ref{l1}.1 and \eqref{e1} imply
\[
 \frac n{\varphi(n)}\sum_{1\le h\le mR}v(h)
 \ll mR\exp\!\left(\sum_{p\mid n}\frac1p
                   +\sum_{p\le mR}\frac{v(p)-1}{p}\right)\ll mR\exp\!\left(\sum_{\substack{p\mid n\text{ with }p>mR}}\frac1p\right)
 =mR e^{L_{mR}(a,b)}.
\]
Finally, it follows from \eqref{e2} and \ref{c1} that 
\[
 \begin{aligned}
 V=\sum_{1\le |h|\le mR}c(h)\le\frac{2\varphi(a)\varphi(b)}m\,
       \frac n{\varphi(n)}\sum_{1\le h\le mR}v(h)\ll R\varphi(a)\varphi(b)e^{L_{mR}(a,b)},
 \end{aligned}
\]
as required.
\end{proof}

We next sum over the denominators to bound the number of ordered pairs
of distinct fractions in $\Fr(\mathcal D)$ at circular distance at most
$\delta/\sum_{a\in\mathcal D}\varphi(a)$. 

\begin{lemma}\label{l3}
	There exists a constant $C_0>0$ such that, for every nonempty
	finite subset $\mathcal D\subset \N_{\ge 2}$ and every
	$0<\delta<1/4$, we have
	\[
	\#\left\{(f_1,f_2)\in \Fr(\mathcal D)^2\colon
	 0<|f_1-f_2|\le \frac{\delta}{\sum_{a\in\mathcal D}\varphi(a)}
	 \right\}
	 \le C_0 \sqrt{\delta}\sum_{a\in\mathcal D}\varphi(a),
	\]
	where $|f_1-f_2|$ denotes circular distance on $\T$.
\end{lemma}

\begin{proof}
	Let $P$ denote the number of ordered pairs in the statement.
	For $a,b\in\mathcal D$, let $V(a,b)$ be the number of ordered
	pairs counted by $P$ whose denominators are $a,b$, respectively. 
	Thus, by summing over the denominators, we obtain
	\[
	 P=\sum_{a,b\in\mathcal D}V(a,b).
	\]
	
    Set $J=\sum_{a\in\mathcal D}\varphi(a)$ and $R=\delta/J<1/4$.
    Since $V(a,b)=0$ when $R\lcm(a,b)<1$, Lemma \ref{l2} implies
    \begin{equation}\label{e6}
     RP\ll\sum_{\substack{a,b\in\mathcal D \text{ with }\\R\lcm(a,b)\ge1}}
       R^2\varphi(a)\varphi(b)e^{L_{R\lcm(a,b)}(a,b)}
     \le e^2\delta^2+
       \sum_{\substack{a,b\in\mathcal D \text{ with }\\R\lcm(a,b)\ge1 \text{ and }\\L_{R\lcm(a,b)}(a,b)>2}}
       R^2\varphi(a)\varphi(b)e^{L_{R\lcm(a,b)}(a,b)}.
    \end{equation}
    
Let $a,b\in \mathcal D$ with $R\lcm(a,b)\ge1$ and $L_{R\lcm(a,b)}(a,b)>2$. It follows from the definition of $L_t(a,b)$ that 
\[
 L_{\lceil R\lcm(a,b)\rceil}(a,b)
 \ge L_{R\lcm(a,b)}(a,b)-\frac12>1.
\]
Indeed, increasing the threshold from $R\lcm(a,b)$ to its ceiling
removes at most one prime from the sum, and its reciprocal is at most $1/2$.
Since $L_t(a,b)$ is nonincreasing in $t$ and eventually zero,
there is a unique integer $j\ge\lceil R\lcm(a,b)\rceil$ such that
$L_j(a,b)\ge1>L_{j+1}(a,b)$, whence, by \eqref{e1},
\[
 L_{R\lcm(a,b)}(a,b)
 =L_{j+1}(a,b)+
 \sum_{\substack{p \text{ divides } ab/\gcd(a,b)^2 \text{ with }\\ R\lcm(a,b)<p\le j+1}}\frac1p\le1+\sum_{h=2}^{j+1}\frac1h
 \le1+\log(j+1).
\]
Let $\mathcal D_j=\{(a,b)\in \mathcal D^2\colon R\lcm(a,b)\ge1,\ L_j(a,b)\ge1>L_{j+1}(a,b), \text{ and } L_{R\lcm(a,b)}(a,b)>2\}$. Then  every $(a,b)\in \mathcal D_j$ satisfies
\begin{center}
$R\lcm(a,b)\le j$,  \quad $L_j(a,b)\ge1$, \quad and\quad  $e^{L_{R\lcm(a,b)}(a,b)}\le e(j+1)$,
\end{center}	
Applying Lemma \ref{l1}.2 with $y_0=j/R$ and $t=j$, we obtain
\begin{equation}\label{e7}
 \begin{aligned}
 \sum_{(a,b)\in\mathcal D_j}
   R^2\varphi(a)\varphi(b)e^{L_{R\lcm(a,b)}(a,b)}
 &\le e(j+1)R^2
   \sum_{\substack{a,b\in\mathcal D\text{ with }R\lcm(a,b)\le j \text{ and }\\  L_j(a,b)\ge1}}
       \varphi(a)\varphi(b)\\
 &\ll(j+1)R^2e^{-j}\sqrt{j/R}\,J^{3/2}
   =(j+1)e^{-j}\sqrt j\,\delta^{3/2}.
 \end{aligned}
\end{equation}
Since $\{(a,b)\in \mathcal D^2\colon R\lcm(a,b)\ge1 \text{ and } L_{R\lcm(a,b)}(a,b)>2\}=\cup_{j=1}^{\infty} \mathcal D_j$, by
summing over $j\ge1$, it follows from \eqref{e6} and \eqref{e7} that
$
	 RP\ll\delta^2+
	 \delta^{3/2}\sum_{j\ge1}(j+1)e^{-j}\sqrt j
	 \ll\delta^{3/2},
$
	since the series converges and $0<\delta<1/4$. Therefore
$
	 P\ll\frac{\delta^{3/2}}R
	 =\sqrt\delta\sum_{a\in\mathcal D}\varphi(a),
$
	as required.
\end{proof}

\section{Proof of Lemma \ref{l4}}\label{s4}

We use Lemma \ref{l3} to find sufficiently many discontinuities with
disjoint intervals on which the required inequality holds.

\begin{proof}[Proof of Lemma \ref{l4}] Let $F_B(x)=\sum_{b\in B} f(bx)$.
	We first describe the discontinuities of $F_B$.
\begin{claim}\label{c3}
There is a nonempty finite set $\mathcal D\subset\N_{\ge2}$ such that
the set of discontinuities of $F_B$ is $\Fr(\mathcal D)$. At each
$x_0\in\Fr(\mathcal D)$, we have
\[
 F_B(x_0+)-F_B(x_0-)\in\mathbb Z\setminus\{0\}.
\]
\end{claim}

\begin{proof}[Proof of \ref{c3}]
For each $b\in B$, the function $x\mapsto f(bx)$ is right-continuous
and locally constant except at the points satisfying
$bx\equiv1/3$ or $2/3\pmod1$. Thus every discontinuity of $F_B$
belongs to the finite set
\[
 \bigcup_{b\in B}
 \left\{\frac{3r+1}{3b},\frac{3r+2}{3b}:r=0,\ldots,b-1\right\}.
\]
In particular, every discontinuity has a rational representative.

For each $x_0\in\T$, write
\[
 F_B(x_0+)=\lim_{\varepsilon\to0^+}F_B(x_0+\varepsilon)
 \quad\text{and}\quad
 F_B(x_0-)=\lim_{\varepsilon\to0^+}F_B(x_0-\varepsilon),
\]
where addition is taken modulo one. Both limits exist. Multiplication
by a positive integer preserves the direction of approach, so
\[
 \di_b(x_0):=\lim_{\varepsilon\to0^+}
 \bigl(f(b(x_0+\varepsilon))-f(b(x_0-\varepsilon))\bigr)
 =\begin{cases}
 1,&\text{if }bx_0\equiv1/3\pmod1,\\
 -1,&\text{if }bx_0\equiv2/3\pmod1,\\
 0,&\text{otherwise}.
 \end{cases}
\]

Let $\chi$ be the real character modulo $3$, with values
$\chi(0)=0$, $\chi(1)=1$, and $\chi(2)=-1$ on the residue classes.
For a reduced fraction $u/d$, we obtain that $\di_b(u/d)\neq 0$ only if $d\mid 3b$.
When $d\mid3b$, the formula above implies
$\di_b(u/d)=\chi(3bu/d)=\chi(u)\chi(3b/d)$, including the case
where this value is zero. Summing over $b\in B$, we obtain
\begin{equation}\label{e9}
 F_B(u/d+)-F_B(u/d-)
 =\chi(u)\sum_{\substack{b\in B\text{ with }d\mid3b}}\chi(3b/d).
\end{equation}
To see that a discontinuity exists, take $m=\max B$ and $d=3m$.
For $b\in B$, the condition $3m\mid3b$ holds only when $b=m$.
Hence \eqref{e9}, with $u=1$, implies
\[
 F_B\left(\frac1{3m}+\right)-F_B\left(\frac1{3m}-\right)
 =\chi(1)\chi(1)=1.
\]
In particular, the point $1/(3m)$ is a discontinuity.

The function $F_B$ is right-continuous, so its value at $u/d$ equals
its right limit. It is therefore discontinuous there precisely when
the difference in \eqref{e9} is nonzero.

If $3\nmid d$, the condition $d\mid3b$ implies $d\mid b$, and hence every
term $\chi(3b/d)$ is zero. If $3\mid d$, the condition $\gcd(u,d)=1$
implies $\chi(u)\in\{-1,1\}$. Therefore
\[
 F_B\text{ is discontinuous at }u/d
 \quad\Longleftrightarrow\quad
 3\mid d\text{ and }
 \sum_{\substack{b\in B\text{ with }d\mid3b}}\chi(3b/d)\ne0,
\]
where the sum is independent of the numerator $u$. Define
\[
 \mathcal D=\left\{d\in\N:3\mid d\text{ and }
 \sum_{\substack{b\in B\text{ with }d\mid3b}}\chi(3b/d)\ne0\right\}.
\]
The set of discontinuities is exactly $\Fr(\mathcal D)$, so
$\mathcal D$ is nonempty. For each $d\in\mathcal D$, the nonzero
sum $\sum\chi(3b/d)$ contains a term indexed by some $b\in B$ with $d\mid3b$.
Thus $d\le3\max B$, and the set $\mathcal D$ is finite. Finally,
the nonzero difference in \eqref{e9} is an integer because
all values of $\chi$ are integers.
\qedhere[End of \ref{c3}]
\end{proof}

	Let $\mathcal D$ be as in \ref{c3}, let
$J=\#\Fr(\mathcal D)=\sum_{a\in\mathcal D}\varphi(a)$,
	and let $R=\delta/J$, where $0<\delta<1/4$ will be chosen below.
	Let $P_R$ count the ordered pairs of distinct points of
	$\Fr(\mathcal D)$ at distance at most $R$. Then Lemma \ref{l3} implies that
	\[
	P_R\le C_0\sqrt\delta\,J.
	\]
	Choose $0<\delta<1/4$ such that $C_0\sqrt\delta\le1/2$.
	Then at most $J/2$ discontinuities of $F_B$ have another discontinuity
	within distance $R$. Thus at least $J/2$ discontinuities have no
	other discontinuity within distance $R$.
	For each of these at least $J/2$ discontinuities $x_0$, it follows from \ref{c3} that
	 $F_B(x_0+)-F_B(x_0-)\in\mathbb Z\setminus\{0\}$, whence
	\begin{center}
	 $|F_B(x_0+)|\ge1/2$\quad  or\quad $|F_B(x_0-)|\ge1/2$.
	 \end{center}
	There is no other discontinuity within  $(x_0-R, x_0+R)$,
	so $F_B$ is constant on each of the intervals $(x_0-R/2, x_0)$ and $(x_0,x_0+R/2)$. On at least one
	of these intervals, we have $|F_B(x)|\ge1/2$.
	Select one such interval for each of these discontinuities.
	The selected intervals are disjoint up to endpoints, since the chosen points
	have pairwise distance greater than $R$. Consequently,
	\[
	\mu\{x\in\T:|F_B(x)|\ge\tfrac12\}
	\ge J/2\cdot  R/2=\delta/4.
	\]
	This proves the lemma with $\rho=\delta/4$.
\end{proof}

\section*{Acknowledgement}
This work was supported by the National Natural Science Foundation of China (Grant No. 12571578) and by Xidian University (Grant No. 21103260003).

\end{document}